\documentclass[12pt]{amsart}
\usepackage{tcolorbox}
\usepackage{multirow}
\usepackage{graphicx}
\usepackage{amssymb}
\usepackage{stmaryrd}
\author{M.  Shevtsova}
\address{\textbf{Margarita Shevtsova:} Department of Discrete Mathematics, Moscow institute of physics and technology,  Dolgoprudny 141700, Russia.}
\email{shevtsova.ma@phystech.edu}
\thanks{}

\author{A. Kanel-Belov}
\address{\textbf{Alexei Kanel-Belov:} Department of Discrete Mathematics, Moscow institute of physics and technology,  Dolgoprudny 141700, Russia.}
\email{kanelster@gmail.com}
\thanks{}

\author{M. Golafshan}
\address{\textbf{Mehdi Golafshan:} Department of Discrete Mathematics, Moscow institute of physics and technology,  Dolgoprudny 141700, Russia.}
\email{m.golafshan@phystech.edu}
\thanks{}

\title{{An indirect method for solving combinatorial problems}}

\date{}
\begin{document}

\begin{abstract}
This article introduces a pedagogical method for {\it solving combinatorial problems} that frequently involve structures that are unfamiliar or less familiar.
Indeed, an indirect method has been proposed in order to evade any possible questions that might obstruct our path to the correct answer and final solution. 
With this method, the desired configuration can be formulated.
This method is effective for solving problems in Olympiads and mathematical competitions as well as for teaching combinations and counting science in the elementary school mathematics curriculum.

\end{abstract}

\maketitle

\noindent\rule{12.7cm}{1.0pt}

\medskip


\noindent\rule{12.7cm}{1.0pt}

\section{Introduction}

Combinatorics is a significant subfield of discrete mathematics, which is defined as “a vibrant branch of contemporary mathematics that is widely applied in business and industry" (National Council of Teachers of Mathematics [NCTM], 2000, p. 31).
Combinatorics provides tools for dealing with both everyday life and professional practice, and it is connected to numerous branches of mathematics and other disciplines (for example, computer science, communication, genetics, and statistics). Thus, it is critical to integrate combinatorics into the mathematics curriculum beginning in the early elementary grades and continuing through senior high school (English, 1993; NCTM, 2000).

\medskip

Although, combinatorial problems from mathematical olympiads require a high level of wit and creativity to solve. There are difficult problems that can be solved almost entirely by a single brilliant idea. In order to approach these problems with confidence, it is necessary to have previously dealt with problems of similar difficulty.
A large number of combinatorial problems have a quick and/or simple solution. However, this does not mean the problem at hand is not difficult. The difficulty of a problem in combinatorics is frequently due to the fact that the solution is well “hidden".
In recent years, good and coherent books have been written in this field, among which we can mention \cite{1}.
Moreover, articles such as \cite{2} have been written on this subject.

\medskip

In addition, solving combinatorial problems is related to the nature of combinatorial thinking.
It is discussed in greater detail in article \cite{3}.

\subsection{ Purpose of the study}

Combinations, or discrete mathematics in general, are a type of science, though some of the problems they solve have ancient mathematical roots. Problem solving and teaching are both critical components of learning this branch of mathematics. We attempt to introduce a method that, to some extent, lifts the veil of secrecy surrounding combinatorial problems.

\medskip

Therefore, this research has two general and serious goals:

\begin{enumerate}
    \item[\textbf{1)}] 
    Developing and implementing a pedagogical and methodical approach to teaching combinatorial concepts and solving combinatorial problems for teachers;
    
    \item[\textbf{2)}]
    Development of methods for solving combinatorial problems in Olympiads.
\end{enumerate}

\subsection{Mathematical backgrounds}
Here, establishes the terminology and concepts that will be used throughout the paper.
Though the majority of you probably already know it,  yet their summary is not without charm.

\medskip

\large
Notations:
\normalsize
We are solely concerned with finite objects. We employ the conventional set-theoretic notation.
Let
$A, B$ are two finite sets.
$|A|$ denotes the {\it size} (the {\it cardinality}) of a set $A$
Moreover,
$A \times B = \{(a, b) : a \in A, b \in B\}$ 
is 
{\it Cartesian product} of $A$ times $B$.
In set-builder notation,
$A  \cup B = \{ x: x \in A \text{  or  } x \in B\}$.

A {\it permutation} of $A$ is a  bijection $f : A \to A$.

\medskip

\large
Terminologies:
\normalsize
We will discuss briefly some concepts in enumerative combinatorics.

Let $A,B$ be finite sets that are
pairwise disjoint. Then
\begin{itemize}
    \item \textbf{Addition principle:}
$|A \cup B| = |A|+|B|.$
\item \textbf{Product principle:}
$|A \times B|= |A| \cdot |B|.$
\end{itemize}

\newpage
We propose a slightly more complicated definition.
Then,
the number of permutations of $n$ objects with $n_1$ identical objects of type $1$, $n_2$ identical objects of type $2$, . . . , and $n_k$ identical objects of type $k$.
\begin{itemize}
    \item \textbf{Permutations with repetition:}
    $\frac{n!}{n_1!n_2! \cdots n_k!}.$
\end{itemize}

Refer to references \cite{4} or \cite{5} for additional information and study.

\subsection{Methodology }
We speculated as mathematics teachers and researchers about a claim that learning combinatorial concepts necessitates a unique way of thinking, and a review of the related literature revealed that some researchers acknowledged this speculation and labeled it combinatorial thinking.
Investigating certain situations, systematically producing all cases, converting the problem into another combinatorial problem, and understanding combinatorial reasoning were the four levels of understanding combinatorial thinking we utilized.

The assessment of relevant literature in the combinatorial problem convinced us that much more research is needed to illuminate the teaching and learning of combinatorics as an engaging aspect of both school and university mathematics curricula.

The authors drew on their previous teaching experiences to conduct their research.
It is a kind of {\it reverse engineering} that we are trying to do here.
We ask the following general question at each step of solving the relevant problem:

\begin{tcolorbox}
    \textbf{{What is the source of my blunder and my confusion in solving this problem?}}
\end{tcolorbox}

\subsection{Strategy}
When we are faced with a problem involving combinations that are not the same as the previous ones we have solved, we employ this strategy and ask the this question. At the end, in each step to prevent mistakes, we will try to design a problem for the answer and see if the answer corresponds to the main problem one by one or not.

\medskip

To improve our proficiency with this approach and to make our discussion more coherent, we begin by examining two problems in this manner. Then, we apply the same strategy to two additional problems and summarize our findings in the discussion section.

\newpage
\section{Finding}

We start by gradually expressing the desired approach through a few examples, and then we look at each one separately.

\medskip

\noindent
\textbf{Sample 1.}
Look at the points in \textsc{Figure} 1.
 How many squares can be formed by connecting points with lines that are parallel to the square's sides? (With the exception of the outer square.)
\begin{figure}[ht]
    \centering
    \includegraphics[width=0.45\textwidth]{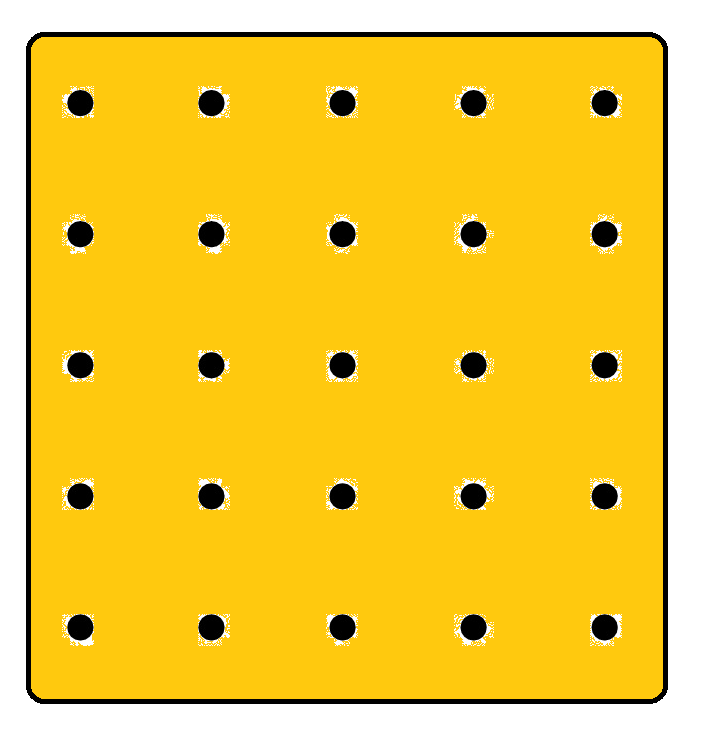}
    \caption{Set of integer points}
\end{figure}

\medskip

Our main goal is not to solve this problem and find an answer; rather, we want to establish a methodical and purposeful approach that we can apply to other problems.

\medskip

\textbf{Step i)}
This can be a kind of wording problem.
In fact, we only have one shape, and we are not sure what to count as part of it.
What we do know is counting the number of squares.
So, the first question that may cause us to become perplexed is:
\footnote{The word configuration is used to describe a finite collection $p_1 \ldots p_n$ of points in $\mathbb{R}^{2}$ and some lines.}
\begin{tcolorbox}
    \textbf{What is the configuration
    of the squares in the figure?}
\end{tcolorbox}

To answer this question, one must act as if it were a wording problem, establishing a logical link between the problem form and a mathematical concept embedded in the figure.

Roughly speaking,
a square is made up of four equal line segments with the facing segments parallel to each other, as we know from basic geometry (of course, considering the inside, which we did not mention here because it has no effect on the answer to the problem).
So regardless of the size of the squares, we are looking for something like \textsc{Figure} 2:
\begin{figure}[ht]
    \centering
    \includegraphics[width=0.50\textwidth]{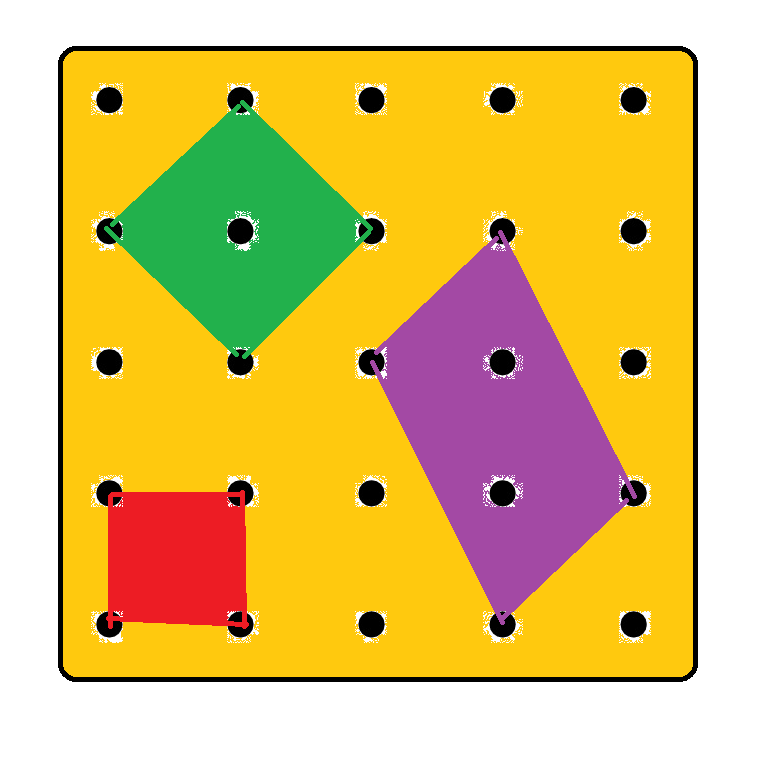}
    \caption{Square in different geometric structure}
\end{figure}

Here are three different colored shapes in purple, red, and green
(And maybe black points!).
Each of these is considered a square in different categories, yet what is meant is a red square (why?)
\footnote{The discussion of the purple square goes back to concepts of Integer Geometry that are beyond the scope of basic mathematics. 
You can look at \cite{6} for further details.}.
Another simple, short, and crucial thing is required in order to more precisely determine your desired configuration.
We must once again determine what caused our mistake and respond accordingly.

\medskip

\textbf{Step ii)}
This section could have been a subsection of the previous section, but because it is so important, we will look at it separately.
What we will do now is specify the configuration type. However, this necessitates a thorough examination of the problem statement.
This is what confuses us:
\begin{tcolorbox}
    \textbf{Which of the following squares should be counted according to the problem (red or green or both)?}
\end{tcolorbox}

\medskip

Let us go over the problem once more, and then bold the section where the distinction between green and red squares is made.

\medskip

The points marked in the following Cartesian coordinate system are integer.
\begin{tcolorbox}
    \textbf{How many squares with
    sides that are parallel to the coordinate axes and vertices that are these points can be built?
    }
\end{tcolorbox}

\medskip

Therefore, due to the bold part, we only have to consider the red squares.
But from a pedagogical point of view, it is good to ask the question, what would the problem be if we were to count the green squares as well?
Ask students to discuss this with each other in class.
Students can learn and consolidate what they have learned by creating a problem based on various data and answers.
But the problem that can be raised here is a dual problem:

\medskip

\noindent
\textbf{Sample $1'$.}
Look at the points in \textsc{Figure} 3.
How many squares with  vertices that are these points can be built?
(
To see the solution to this problem, you can refer to  \cite{7}.)

\begin{figure}[ht]
    \centering
    \includegraphics[width=0.45\textwidth]{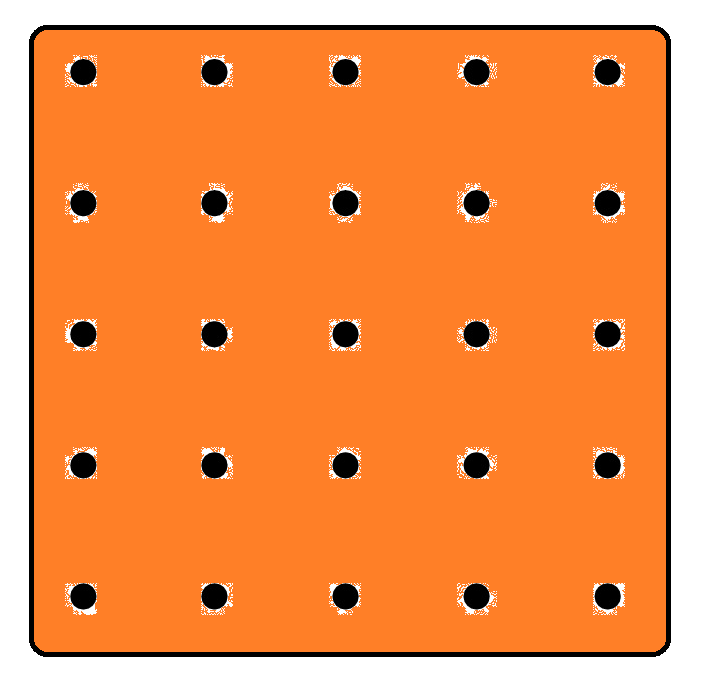}
    \caption{Set of integer points}
\end{figure}

So far, we have discovered that regardless of the square size, we are looking for the  configuration
in \textsc{Figure }4.

\begin{figure}[ht]
    \centering
    \includegraphics[width=0.18\textwidth]{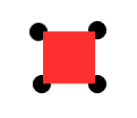}
    \caption{The configuration structure that is considered.}
\end{figure}

\medskip

\newpage
Now that we have a clear understanding of the configuration structure, regardless of the size of the squares, we can move on to this.

\medskip

\textbf{Step iii)}
In this section, we will review the allowed sizes of squares.
We now have a very large setup of objects and conditions. So, in order to make counting easier, we must break it down into smaller structures.
These smaller structures are squares and can come in a variety of sizes.
However, we are once again perplexed by the following question:
\begin{tcolorbox}
    \textbf{ What is the allowable size for squares?
    }
\end{tcolorbox}

\medskip

The simplest criterion that determines the size of a square is its side.
 As a result, we must respond to the following question: What are the allowed segments that can be sides of this square?
 If we call each row and each column a line, then we can say that the points on all the lines are equal.
 As a result, the type of square is determined by the line segment that is the square's floor side(
 See \textsc{Figure} 5).

\begin{figure}[ht]
    \centering
    \includegraphics[width=0.70\textwidth]{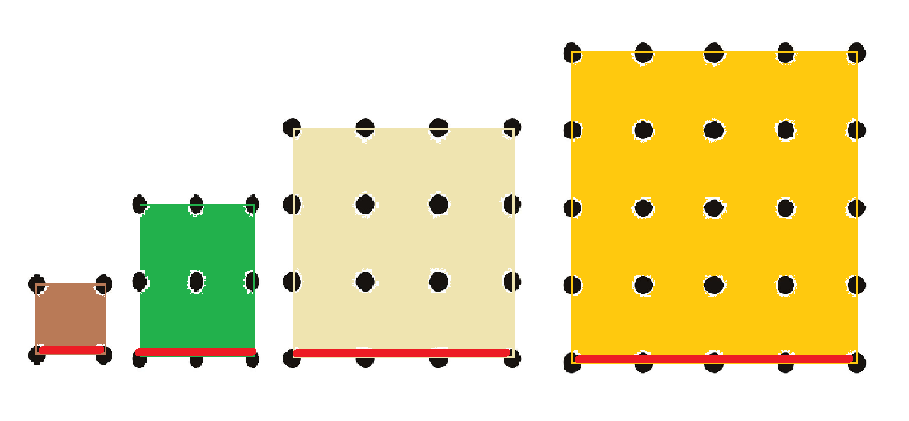}
    \caption{Different sizes of squares in the configuration}
\end{figure}

\medskip

According to the above and the figure, we find that the minimum distance between two different horizontal points that can be the sides of a square is equal to $1$, and the maximum distance between two horizontal points that can be the sides of a square is equal to $4$.
The desired configuration is now fully revealed to us after three steps.

\medskip

In the final discussion, fourth, we count the final number of squares.

\medskip

\textbf{Step iv)}
First, for the convenience of work, we show a square with sides $1, 2, 3$ and $4$ with the symbols of $A_1, A_2, A_3$ and $A_4$, respectively.
So the question in mathematical form is as follows:
Find the $|A_i|$ for all $i=1,2,3,4$.
What is the fundamental and ambiguous question that arises for each of these sets?
After counting, we ask ourselves whether we have counted too much or too little?
Therefore, the following puzzling question arises for us:
\begin{tcolorbox}
    \textbf{ How do we make sure that it is true? 
    }
\end{tcolorbox}

\medskip

If we are not experts at solving combinatorial problems and are unfamiliar with the process, we may count each square in $A_i$ individually.
As an example, let's say we did this for $A_2$.
In this process, the number of squares of size $2$ may be more or less than counted.
But how do we know if this number is less than or greater than the correct number?
Attempting to create a problem for the answer we provide is a good guideline.
For example, if you found $8$ squares of size $2$, ask the question: If we had $8$ squares of size $2$, would this be our configuration?
We will find a one-to-one correspondence between the answer and the desired problem if we proceed in this manner.
However, for the sake of convenience and to avoid mistakes, we can ask the following questions:
\begin{tcolorbox}
    \textbf{ How do you begin counting and how do you count according to specific rules?
    }
\end{tcolorbox}

\medskip

To do this, we do the following: to count the number of elements $A_2$: In the figure, consider virtual horizontal rails at a distance of two units
( See \textsc{Figure} 6).

\begin{figure}[ht]
    \centering
    \includegraphics[width=0.9\textwidth]{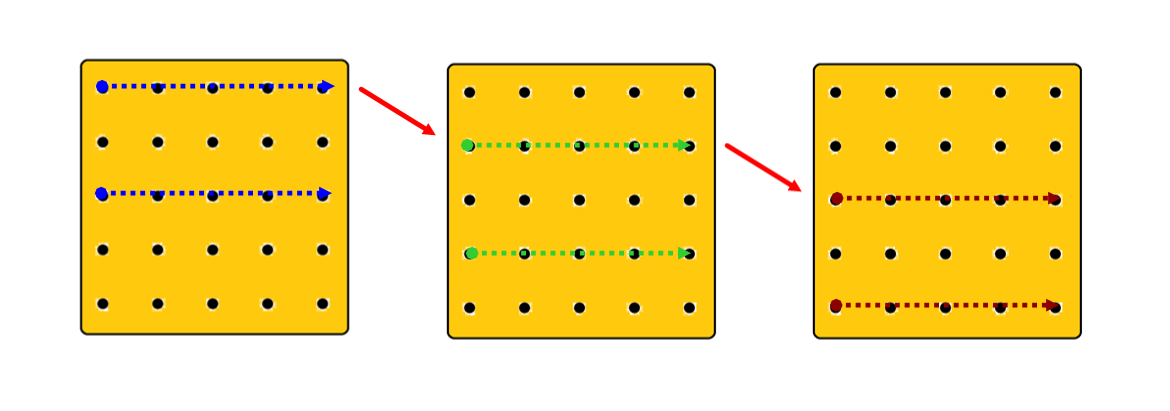}
    \caption{Structures of  squares with a side of two units}
\end{figure}

\medskip

According to the figures from left to right, if we move the rails one unit down, we have $3$ pairs of these virtual rails.
The number of squares between the blue, green, and brown rails is indicated by the symbols $A_{21}, A_{22}$ and $A_{23}$, respectively.
Hence,
$|A_2| = |A_{21}|+|A_{22}|+|A_{23}|.$
That is, we use the principle of addition. Due to the symmetry of the shape in terms of the number of points as well as the rails, we count only the number of squares created by the blue rail ( See
\textsc{Figure} 7).

\begin{figure}[ht]
    \centering
    \includegraphics[width=0.85\textwidth]{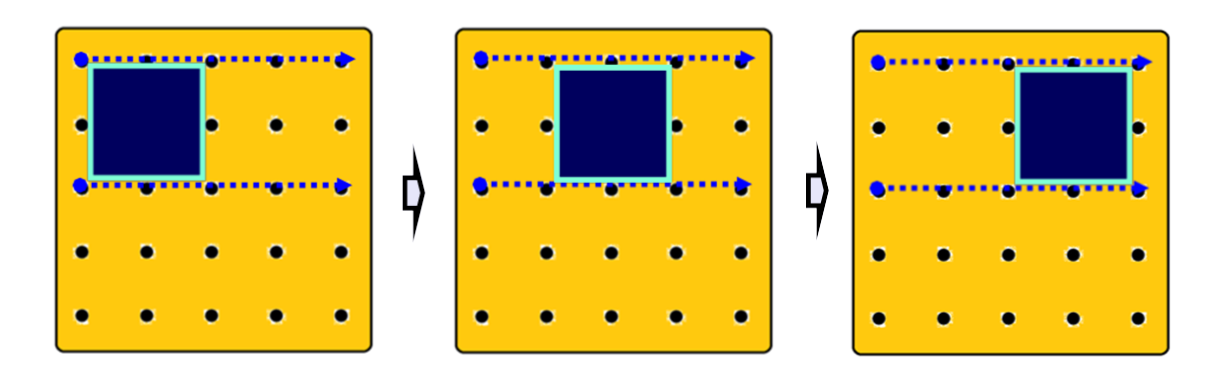}
    \caption{Counting squares with sides of two units in the first structure}
\end{figure}

According to the figure,
$|A_{21}|=3$,
and so on
$|A_{22}|=|A_{23}|$.
Then
$|A_2|=9$.
If we repeat the process for virtual rails at distances of $1$, $3$, and $4$, we get the following results.
The number of pairs of rails at a distance of one is $4$, and between each of these rails we have $4$ squares on one side; So according to the principle of multiplication we have:
$|A_1| = 4 \cdot 4$.
Also,
$|A_3|= 2 \cdot 2$, and
$|A_4| = 1 \cdot 1$.

Finally, based on the principle of addition, the total number of squares equals:
$|A|=|A_1|+|A_2|+|A_3|+|A_4|$
i.e.,
$|A|=30$.

\medskip

To help you better understand and practice this method, we will look at another combinatorial problem together.
Naturally, we chose a less well-known problem this time.
As a result, the following problem will be examined in this manner:

\medskip

\noindent
\textbf{Sample 2.}
What is the number of all possible ways to read ``\textbf{Open!}'' in \textsc{Figure} 8, according to these conditions:
\begin{enumerate}
    \item[\textbf{a)}] the square with next letter must share a side with the previous square;
    \item[\textbf{b)}] the sequence of letters in the read squares must make up the word ``\textbf{Open!}''.
\end{enumerate}

\begin{figure}[ht]
    \centering
    \includegraphics[width=0.28\textwidth]{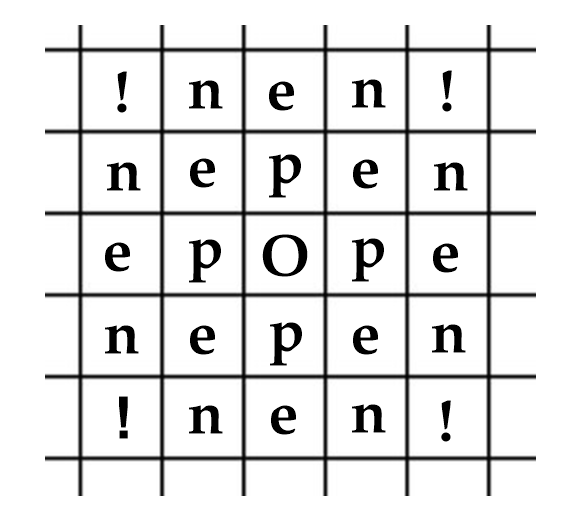}
    \caption{Word game table
}
\end{figure}

\newpage
We attempt to replicate the method used to solve the previous problem and reach a general conclusion about this one in particular.

\medskip

\textbf{Step i)}
This is an intriguing problem that also serves as an excellent educational example. While the symmetrical structure of this configuration may be perplexing at first, it ultimately aids in our solution.
In this problem, we have been asked to come up with alternative ways to read the word ``\textbf{Open!}'' in the hypothetical form given the above conditions.
The first step is to establish a meaningful and palatable link between the form and the problem. As a result, the first question that may cause us to become perplexed is:
\begin{tcolorbox}
    \textbf{ What is the configuration of the word ``\textsc{Open!}'' in the figure?
    }
\end{tcolorbox}

\medskip

We first examine the word structure, then do the same for the shape, and finally make a connection between them.

This word contains five symbols, and more importantly, it begins with the ``\textbf{O}'' symbol and concludes with the ``\textbf{!}'' symbol
(See \textsc{Figure } 9).
\begin{figure}[ht]
    \centering
    \includegraphics[width=0.28\textwidth]{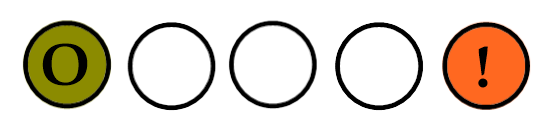}
    \caption{Main configuration}
\end{figure}

A circle indicates the location of the symbols. For the time being, we can say that this structure exists independently of the figure's main configuration, and that the general structure reveals its nature when it interacts with the shape.

In \textsc{Figure} 10, we attempt to specify one or more of the possible readings of the word ``\textbf{Open!}'' by specifying the initial and final letters:
\begin{figure}[ht]
    \centering
    \includegraphics[width=0.62\textwidth]{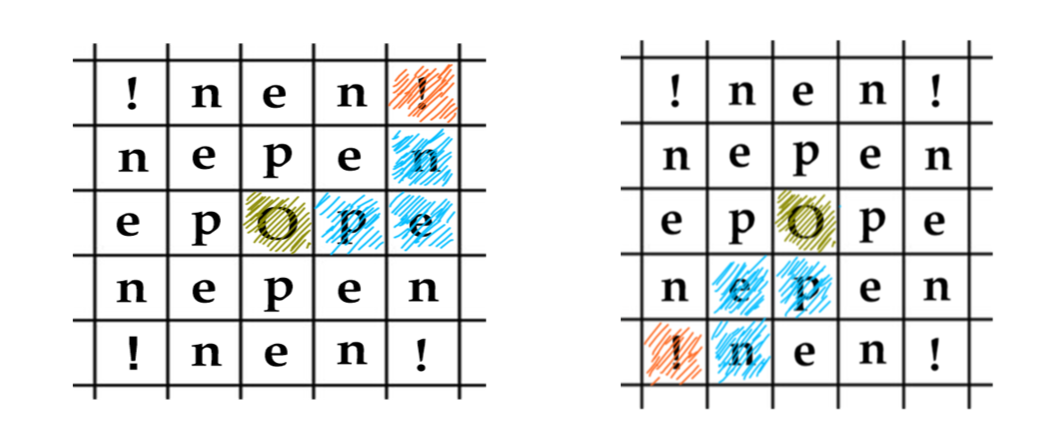}
    \caption{Types of configuration}
\end{figure}

\newpage

Thus far, the word has been configured in part. 

\medskip

\textbf{Step ii)}
What may cause confusion here (though hopefully not too much!) is whether this is the correct configuration in the form that corresponds to the word in question. Thus, the issue that concerns us here is as follows:
\begin{tcolorbox}
    \textbf{Which cells can be configured in the desired manner and which cells cannot?
    }
\end{tcolorbox}

\medskip

As previously stated, it is critical to pay close attention to the problem's assumptions and conditions. The cells whose adjacent symbols should be adjacent are mentioned in this case. What may cause students to become confused in this instance is the term ``adjacent."
Adjacent here, the existing cells of two consecutive symbols must have a common side. That is, structures as follows are not allowed
(
\textsc{Figure}
11).

\begin{figure}[ht]
    \centering
    \includegraphics[width=0.30\textwidth]{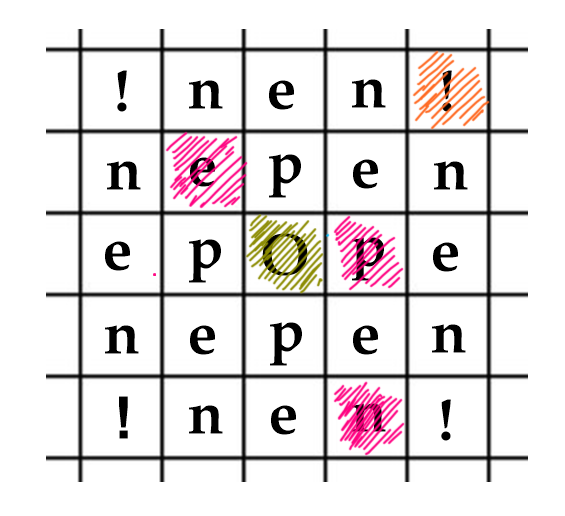}
    \caption{Unacceptable configurations}
\end{figure}

Solicit from students the creation of a problem for this structure that is unrelated to our configuration.
One of the problems that can be raised is the following:

\medskip

\noindent
\textbf{Sample $2'$.}
What is the number of all possible ways to read ``\textbf{Open!}'' in the picture such that,
 the sequence of letters in the read squares must make up the word ``\textbf{Open!}''.

\medskip

The configuration has been precisely defined. That is, it is critical that the signs are adjacent conditionally and have no bearing on their proximity from the left, top, or right. That is, the word can be read both left to right and right to left.

\medskip

\textbf{Step iii)}
We organize the authorized structures in this section according to the legal system in order to familiarize ourselves with all the configuration structures.
Thus, the question that prompts our skepticism in this instance is as follows:
\begin{tcolorbox}
    \textbf{ What general structure should these cells have in order to represent this word?
    }
\end{tcolorbox}

\medskip

If we pay attention to the following figure, there are only ``\textbf{!}'' symbols in the four corners, and we have only one ``\textbf{O}'' in the center. So the {center square} must be present in our structures
(See
\textsc{Figure} 12).

\begin{figure}[ht]
    \centering
    \includegraphics[width=0.30\textwidth]{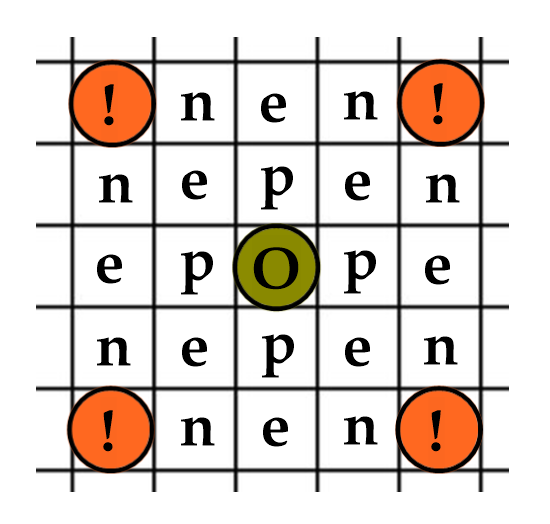}
    \caption{Common configuration features}
\end{figure}

\medskip

This means that the olive cell must appear as the first symbol in each of the configurations, and one of the four orange cells must also appear in our structure.
Now, let us count.

\medskip

\textbf{Step iv)} However, this question perplexes us:
\begin{tcolorbox}
    \textbf{ How and where to begin counting?
    }
\end{tcolorbox}

\medskip

We begin by selecting ``\textbf{O}''.
We fix it because it is required in all structures.
Additionally, the selection of ``\textbf{!}'''s is straightforward. We will use one of them as an illustration.
There are four ways we can do this.

If we number the ``\textbf{!}'''s in this manner 
``\textbf{$\text{!}_1$}'',
``\textbf{$\text{!}_2$}'',
``\textbf{$\text{!}_3$}'',
and
``\textbf{$\text{!}_4$}'',
 the $A_k$ collection is a set of ``\textbf{Open!}'' word configurations in which ``\textbf{$\text{!}_k$}'' is the final symbol.
So if we show the total number of allowed words with $|A|$, according to addition principle we have
$$|A|= |A_1|+|A_2|+|A_3|+|A_4|.$$

\newpage
On the other hand, we consider only one of the structures due to the symmetry of the symbol distribution. They are denoted by a dash in \textsc{Figure} 13.
\begin{figure}[ht]
    \centering
    \includegraphics[width=0.65\textwidth]{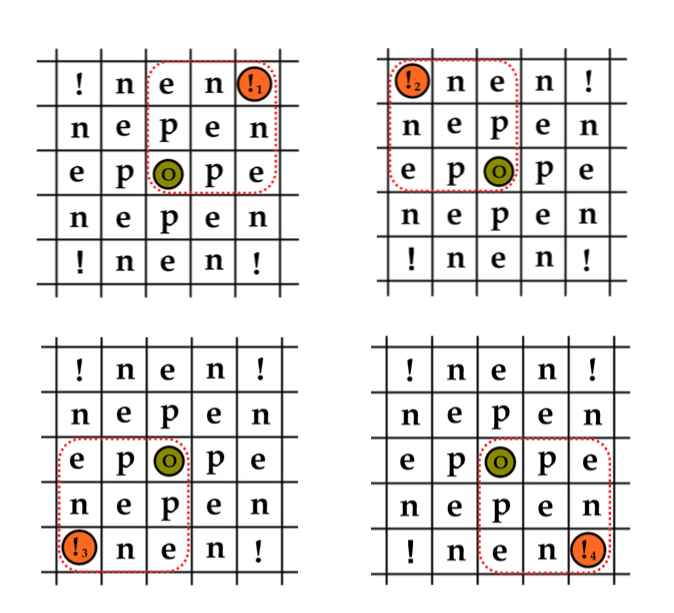}
    \caption{Components of configurations}
\end{figure}

As a result of the multiplication and symmetry principles, we have:
$|A| = 4 \cdot |A_1|.$
As a result, we only need to count one of them. For instance, you could do this for $A_1$.
However, how is $|A_1|$ to be calculated?
According to the figure, in order to get from ``\textbf{O}'' to ``\textbf{!}'' using the desired structure, we must ascend two times and turn right two times. Of course, the order in which the movements occur is irrelevant (why?).

If we depict the upward movement with the $\bf{U}$ and the right movement with the $\bf{R}$, we must count the structure of the $\bf{UURR}$, which equals the repetition:
$|A_1|= {4!}/{2!2!}$,
i.e.,
$|A_1|=6$.
Thus, the total number equals:
$|A| = 4 \cdot 6 $,
i.e.,
$|A|=24$.


\medskip

\section{Discussion}
The following questions could arise if we try to solve the problem using the two examples given:

\medskip

\begin{enumerate}
    \item[\textbf{1)}]
    Identify the objects that need to be counted and look at the configuration of the objects;
      \item[\textbf{2)}]
      Indicate the configuration's attributes;
       \item[\textbf{3)}]
       Look for a structure in the set up (for example, divide objects into classes);
       
       \item[\textbf{4)}]
       Try to use the product or addition principle, and then verify your conclusion by establishing a one-to-one relationship between the items you enumerated and the configuration's component pieces.
\end{enumerate}
 
 \medskip

 We produce a response for each step in order to try to make a problem out of the bad answers and figure out why they were wrong.
 That's why we call it an indirect method ( We have placed these questions in grey boxes.).

 \medskip

 Finally, we try to go over each step once with the following query to make sure we have covered everything:
 
 \begin{tcolorbox}
    \textbf{ How do we ensure that the counting is carried out properly?
    }
\end{tcolorbox}
 
 \medskip

 The steps will undoubtedly be time-consuming and unpleasant, but they will be helpful. Besides, after we become professionals, we are not always required to record every action we conduct.

\section{Conclusion Remarks}
Combinatorics is recognised as one of the more difficult mathematical topics to teach and learn, as we have previously discussed and encountered. Since most issues lack simple fixes, there is a great deal of confusion over how to approach them and which method to employ. There are many instances where two different approaches that result in two different answers to the same problem can seem equally compelling.
While there are other ways to solve combinatorial puzzles besides the one described in this article, it might be helpful when these puzzles are challenging.
At the start of teachers' combinatorics instruction, it can also be employed as a teaching strategy.

Finally, the authors want to investigate more issues relating to teaching discrete and combinatorial mathematics in the future.

\subsection*{Acknowledgement}
 Alexei Kanel-Belov is supported by the Russian Science Foundation 22-11-00177.


\newpage
\begin{thebibliography}{99}

\bibitem{1}
P. Soberón.
{\it
Problem-Solving Methods in Combinatorics-An Approach to Olympiad Problems.}
Birkhäuser (2013).
\textcolor{red}{DOI: 10.1007/978-3-0348-0597-1}.

\medskip

\bibitem{2}
M. M. Eizenberg and O. Zaslavsky.
{\it Students' Verification Strategies for Combinatorial Problems.}
Math. Think. Learn.,
vol.
\textbf{6}
(2004): 15-36.
\textcolor{red}{DOI: 10.1207/s15327833mtl0601-2}.


\medskip

\bibitem{3}
Y. M. Hidayati, C. Sa’dijah and S. A. Qohar.
{\it Combinatorial Thinking to Solve the Problems of Combinatorics in Selection Type.}
Int. J. Learn. Teach. Educ. Res.,
vol.
{\bf 18}
(2019):
65-75.
\textcolor{red}{DOI: 10.26803/IJLTER.18.2.5}.


\medskip

\bibitem{4}
M. Bona.
{\it 
A Walk Through Combinatorics: An Introduction to Enumeration and Graph Theory.
}
World Scientific Publishing
(2006).
\textcolor{red}{DOI: 10.1142/10258}.

\medskip

\bibitem{5}
M. Keller and
 W. Trotter.
 {\it
  Applied Combinatorics.
 }
American Institute of Mathematics
(2017).
\textcolor{red}{ISBN: 13: 9781973702719}.

\medskip

\bibitem{6}
O. N. Karpenkov. 
{\it 
Geometry of Continued Fractions.
}
AACIM,
vol.
{\bf 26},
 Springer,
 Berlin, Heidelberg
 (2022).
\textcolor{red}{DOI: 10.1007/978-3-662-65277-0}.

\medskip

\bibitem{7}
T. Andreescu and
Z. Feng.
{\it
A Path to Combinatorics
for Undergraduates-Counting Strategies.
}
Birkhäuser
(2004).
\textcolor{red}{DOI: 10.1007/978-0-8176-8154-8}.

\end{thebibliography}
\end{document}